\documentclass{article}

\usepackage{amsfonts,amsbsy,amsmath,graphicx}
\usepackage{latexsym}
\pdfoutput=1

\def\intr{\rm int}
\def\tr{\rm trap}
\def\xv{\vec{x}}

\numberwithin{equation}{section}

\begin{document}



\title{Apriori Estimates for
Many-Body Hamiltonian Evolution of Interacting Boson System}

\author{Manoussos G. Grillakis$^{1}$ and Dionisios Margetis$^{1,2}$\\ \\
${}^1${\it Department of Mathematics, University of Maryland}\\ 
{\it College Park, Maryland 20742, USA}\\ {\it mng@math.umd.edu}\\ \\
${}^2${\it Institute for Physical Science and Technology, University of Maryland}\\
{\it College Park, Maryland 20742, USA}\\
{\it dio@math.umd.edu}} 
\date{}
\maketitle


\begin{abstract}
We study the evolution of a many-particle system whose wave function obeys
the $N$-body Schr\"odinger equation under Bose symmetry. 
The system Hamiltonian describes pairwise particle interactions in the absence of an external potential.
We derive apriori dispersive estimates that express the overall repulsive
nature of the particle interactions. These estimates hold for a wide class
of two-body interaction potentials which are independent of the particle number, $N$.
We discuss applications of these estimates to the BBGKY hierarchy
for reduced density matrices analyzed by Elgart, Erd\"os, Schlein and Yau.
\end{abstract}
\vskip5pt



\section{Introduction}

Quantum mechanics dictates that the dynamics of $N$ interacting particles at zero temperature 
be described by a wave function satisfying the time-dependent $N$-body Schr\"odinger equation.
When $N$ is large this description becomes impractical. A common remedy is to reduce the system 
degrees of freedom by replacing the many-body Schr\"odinger equation by effective partial differential 
equations for single-particle wave functions.
Thus, the many-body interaction is approximately replaced by sums of effective one-particle potentials. This approach amounts to a ``mean-field" approximation, 
which is expected to be exact as $N\to \infty$ in an appropriate sense.
This formulation yields nonlinear evolution laws in $3+1$ dimensions which often lead to successful predictions. In particular, the dynamics
of particles with integer spin (Bosons) has been described by nonlinear Schr\"odinger-type equations~\cite{gross61,gross63,pitaevskii61,wu61,wu98}.

The relation of mean-field approaches to the many-particle
Hamiltonian evolution is an area of active research. In settings with Boson symmetry, fundamental questions
concern the rigorous derivation of nonlinear Schr\"odinger-type equations from $N$-body evolution laws 
and the formulation of theories that transcend mean field~\cite{ch1}. Recently, Elgart, Erd\"os, Schlein and Yau 
(EESY)~\cite{esy1,esy2,esy3,esy4} derived rigorously mean-field limits for Bosons
on the basis of BBKGY-type hierarchies for reduced density matrices~\cite{spohn1,spohn2}.

In this paper we study the $N$-particle Hamiltonian
evolution of a Boson system for large yet {\it finite} $N$. We derive certain
apriori estimates for hydrodynamic quantities, and examine some applications of these estimates 
to BBGKY hierarchies formulated for dilute atomic gases~\cite{esy1,esy2,esy3,esy4}. Our goal is
to understand how many-body properties of the Boson system are connected to the mean-field limit.\looseness=-1

The model equation for the particle system reads
\begin{equation}
-i\partial_{t}\psi =H_{N}\psi~, \label{evolution}
\end{equation}
where $\psi$ is the $N$-body wave function and $H_N$ is the Hamiltonian operator. 
This $H_N$ has the form
\begin{equation}
H_{N}=\sum_{a}H_{a} + \sum_{a\not= b}H_{ab}~,\qquad a,\,b=1,\,\ldots\,,N~, \label{Hamiltonian}
\end{equation}
where $a$ and $b$ label the particles, $H_{a}$ is the part of the Hamiltonian acting solely on particle $a$, 
and $H_{ab}$ denotes the (pairwise) interaction between the particles $a$ and $b$. For identical particles,
$H_{a}$ is described by
\begin{equation}
H_{a}=-\Delta_{a}+V_{\tr}(\xv_{a})
\label{potentials-I}
\end{equation}
and the interaction $H_{ab}$ has the form
\begin{equation}
H_{ab}:=\textstyle{\frac{1}{2}}V_{\intr}\big(\xv_{a}-\xv_{b}\big)~,
\label{potentials-II}
\end{equation}
where $\xv_a$ denotes the vector position of particle $a$ ($\xv_a\in\mathbb{R}^3$),
$V_{\tr}$ is a trapping potential which confines the particles in space, and 
$V_{\intr}$ is a two-body potential. 

In our formulation,~(\ref{evolution}) is converted exactly to a set of
hydrodynamic equations. Apriori estimates for finite $N$
are derived here directly from these hydrodynamic equations. Note that by~(\ref{evolution}) we have $\psi(t)=U(t)\psi(0)$
where $U(t)=e^{itH_{N}}$, a unitary operator. For large $N$, this representation of $\psi(t)$
is not particularly useful for apriori estimates and testable predictions.\looseness=-1 

We proceed to discuss elements of the model. The trapping potential, $V_{\tr}$, is
introduced in experiments in order to keep the particles together. 
To avoid unecessary complications such as the effect of bound states in our analysis, we set $V_{\tr}$ equal to zero, 
$V_{\tr}\equiv 0$. The case with a nonzero external potential, $V_{\tr}\not\equiv 0$, 
is the subject of future work.

Next, we comment on the particle interactions, which are crucial for the system evolution.
In our model, the pairwise interaction, $V_{\intr}$, is considered as a known short-range potential.
For dilute gases, which are the primary focus of this investigation,
the particle collisions are sufficiently rare and weak; hence, the distance $|\xv_a-\xv_b|$ 
($a\neq b$) remains mostly large in an appropriate sense. Consequently, 
each particle is influenced only by gross features of $V_{\intr}$. The specific
details of this potential should be immaterial for our purposes. In this vein, by considering particles that 
repel each other, it suffices to assume that
\begin{equation}
H_{ab}=V(|\xv_a-\xv_b|)~;\quad V\ge -c,\ V^{\prime}(s)\leq 0\quad \forall\,s\in \mathbb{R}^+~.\label{eq:repulsion}
\end{equation} 
Note that we consider only two-body interactions, neglecting effects of three-body and
higher-order collisions~\cite{wu58}. 

A key ingredient of the particle model is the scattering length $l$ of the interaction potential, $V_{\rm int}$. This
$l$ encapsulates the gross features of the two-particle scattering~\cite{esy2}. 
Let the zero-energy scattering solution for the class of potentials
$V_{\intr}(\xv)$ that decay sufficiently fast at infinity have the form $1-w_0(\xv)$ where $\lim_{|\xv|\to\infty}w_0=0$;
then, $l$ is defined by $l:=\lim_{|\xv|\to\infty}(|\xv|w_0)$.
Accordingly, $V_{\intr}$ can be expressed in the scaled form~\cite{esy2,lsy,lsy2}
\begin{equation}
V_{\intr}(\xv)=l^{-2}V_{1}\big(\vert\xv\vert/l\big)~,
\label{scattering length}
\end{equation}
where $V_{1}$ is a fixed potential with scattering length equal to unity. 
Thus, the Hamiltonian $H_N$ has two parameters, namely, the number of
particles, $N$, and the scattering length, $l$. 
In the EESY formulation~\cite{esy1,esy2,esy3,esy4}, it is essential to allow for $N\to\infty$ by keeping
$Nl$ fixed. A parameter that enters the mean-field limit of~(\ref{evolution}) is~\cite{esy2}
\begin{equation}
g:= 4\pi Nl~.\label{GP-pmt}
\end{equation}
Thus, setting $g=O(1)$ implies that the scattering length $l$ scales with $N^{-1}$.

The recent work by EESY~\cite{esy1,esy2,esy3,esy4} deserves special attention, because it provides a rigorous
justification of the mean-field limit of~(\ref{evolution}) for Bosons via an averaging procedure. 
In Sec.~\ref{sec:BBGKY} we outline their approach and derive certain apriori estimates. The uniqueness of the solution to the BBGKY hierarchy
in the mean-field limit is also proved via a different argument by Klainerman
and Machedon~\cite{km}. 

A crucial observation within the EESY framework
is that, in order to be able to pass rigorously to the limit as $N$ increases, 
the two-body potential $V_{\intr}$ has to be appropriately scaled with $N$. For example, EESY~\cite{esy1,esy2,esy3,esy4}
posit that $V_{\intr}$ approaches a delta function as $N\to \infty$, which comes from~(\ref{scattering length})
with $l=N^{-1}$. These interactions cannot be treated as a perturbation since they have
an appreciable effect for sufficiently large $N$. We note in passing that in a more elaborate
model $V_{\intr}$ is replaced by the Fermi pseudopotential~\cite{blatt,huangyang,huangetal,lhy,ly,wu58},
an operator that consists of a delta function times a space derivative acting on $\psi$. 
The use of this pseudopotential allows for interesting extensions but causes mathematical difficulties and lies beyond our present scope.\looseness=-1

A few remarks on a mean-field limit of~(\ref{evolution}) and its variants are in order. 
This limit can be derived heuristically under~(\ref{scattering length}) with $l=O(N^{-1})$
by approximately replacing the wave function $\psi$ by the tensor product
\begin{equation}
\psi \approx \prod_{a=1}^N\phi_{a} \label{tensor1}
\end{equation}
where $\phi_{a}:=\phi(t,\xv_{a})$ and $\int {\rm d}\xv\ \vert\phi\vert^2=1$. By manipulation 
of~(\ref{evolution}) and neglect of terms inconsistent with~(\ref{tensor1}),
a cubic nonlinear Schr\"odinger equation is recovered for $\phi(t,\xv)$. This equation reads~\cite{mng}
\begin{equation}
-i\partial_{t}\phi =H\phi +g\vert\phi\vert^{2}\phi -\mu(t)\phi~,\label{GP-eq}
\end{equation}
where $H=-\Delta +V_{\tr}$ and $\mu(t)=(g/2)\int {\rm d}\xv\,\vert\phi\vert^{4}$; cf.~\cite{gross61,gross63,pitaevskii61,wu61,wu98}. 
This form of $\mu$ does not appear in the results of Gross~\cite{gross61,gross63}, Pitaevskii~\cite{pitaevskii61} and EESY~\cite{esy1,esy2,esy3,esy4}
but is in agreement with the derivation by Wu~\cite{wu61,wu98}. This $\mu(t)$ is not essential since it can be absorbed 
into a phase factor for $\phi$.

It is worthwhile mentioning the complementary view on particle dynamics that invokes operator-valued distributions~\cite{berezin,wu61,wu98};
for the periodic case, see e.g.~\cite{huangyang,lhy,ly}. In this context, one defines the Boson
annihilation operator $\Psi(\xv)$ and its adjoint, $\Psi^\ast(\xv)$, for a particle at point $\xv$. Accordingly,
$H_N$ with the Fermi pseudopotential is written in terms of $\Psi$ and $\Psi^\ast$.
The annihilation operator $a_0(t)$ for the Bose condensate, the macroscopic quantum
state in Bose-Einstein condensation, is introduced as the average $a_0(t)=N^{-1/2}\int {\rm d}\vec{x}\ \{\phi^{\ast}\Psi\}$,
where $\phi(t,x)$ is the (one-particle) condensate wave function. The evolution equation for $\phi$ can
be derived by linearization of $H_N$ in $\Psi-a_0$, which is considered as small in some sense,
and enforcement of~(\ref{evolution}) under~(\ref{tensor1})~\cite{dm,wu61,wu98}.
This approach offers physical insight and is amenable to extensions, particularly the inclusion of
higher-order scattering processes~\cite{ch1,lhy,ly,wu58,wu61}.
However, this formalism is less amenable to rigorous treatment. In particular, the introduction of the Fermi
pseudopotential, which is a type of singular distribution, renders the analysis especially difficult. 

Of course, tensor products such as~(\ref{tensor1}) are not approximate solutions in any classical sense. So, an issue
is how to make this type of approximation meaningful. The BBGKY
hierarchy invoked by EESY~\cite{esy1,esy2,esy3,esy4}, for which the Gross-Pitaevskii equation provides a closure condition, 
offers a meaningful scheme for justifying~(\ref{tensor1}): if the Gross-Pitaevskii
equation is satisfied then all the equations in the hierarchy are satisfied as $N\to\infty$. 

To understand the system evolution by connecting macroscopic
variables such as $\phi$ with microscopic quantities for finite $N$, we resort to apriori
dispersive estimates which express the particle repulsions. These estimates are scale 
independent: we derive them for interaction potentials of the form
$V\big(\vert\vec{x}_{a}-\vec{x}_{b}\vert\big)$, omitting the scaling which can be 
inserted in the end. Our basic underlying assumption is given by~(\ref{eq:repulsion}).

Our apriori estimates are useful for the passage to the limit $N\to\infty$, but
also offer insight into the case with finite $N$. The techniques employed in the present work 
originated in~\cite{cgt,ls} in connection to the nonlinear Schr\"odinger equation; see also
~\cite{ckstt,tt}. A new element of the present problem is the presence of a large number of interaction
potentials and the fact that the hydrodynamic momenta 
are {\it not conserved} by the flow. Our analysis addresses the question
in what sense the potentials are repulsive and what are the precise
implications of this repulsive nature. 

The static version of~(\ref{evolution}) is not addressed here. This case was studied by 
Dyson~\cite{dys} and by Lee, Huang and Yang~\cite{lhy}. A mathematical proof of the Bose-Einstein condensation 
for the time-independent case was provided recently by Lieb et al.~\cite{liebseir2,lieb-book,lsy,lsy2}. 

The remainder of the paper is organized as follows. In Sec.~\ref{sec:collapse} we derive an apriori
space-time estimate for the particle density function after the collapse (identification) 
of two particle positions. Boson symmetry is not required for the main result of 
Sec.~\ref{sec:collapse}. In Sec.~\ref{sec:commutator} we restrict attention to Bosons and
derive another estimate for the square of an appropriately averaged (reduced) particle density.
In Sec.~\ref{sec:BBGKY} we outline the EESY approach based on the BGGKY
hierarchy for reduced density matrices~\cite{esy1,esy2,esy3,esy4,km}; and derive apriori estimates in the context
of this formulation. The (Einstein) summation convention for repeated indices is employed throughout this paper
unless it is noted otherwise. 

%
\section{First Estimate by Coordinate Collapse}
\label{sec:collapse}
%

In this section we derive an apriori space-time estimate for the particle density 
of the $N$-body system described by~(\ref{evolution})--(\ref{potentials-II}) when $N$ is large and finite. 
Boson symmetry is not required for our main result but simplifies it considerably, as shown below. 
Our motivation is to investigate in what sense the wave function $\psi$ 
disperses with time. The methodology followed here originates from the work of Lin and Strauss~\cite{ls}.

Specifically, we show that
\begin{equation}
\sum_{a\neq b}\int {\rm d}t\,{\rm d}\vec{X}_{ab}\ \rho(t,\vec{X}_{ab})\le N^2\,\Vert\psi\Vert_{H^{1}}\,\Vert\psi\Vert_{L^2}~,
\label{estimate-I:NB}
\end{equation}
where the particle density function $\rho$ is defined by
\begin{equation}
\rho(t,\xv_{1},\ldots \xv_{N}):={1\over 2}
\big\vert\psi(t,\xv_{1},\ldots \xv_{N})\big\vert^{2}~.\label{density}
\end{equation}
The reduced vector $\vec{X}_{ab}$ entering~(\ref{estimate-I:NB}) 
comes from a basic collapse mechanism which identifies two of the
variables $(\vec{x}_{1},\ldots \vec{x}_{N})$. Hence, $\vec{X}_{ab}$ is defined by
\begin{equation}
\vec{X}_{ab}:=\big(\vec{X}\big)_{\vec{x}_{a}=\vec{x}_{b}}
\in\mathbb{R}^{3(N-1)}~,\qquad \vec{X}:=(\vec{x}_1,\,\vec{x}_2,\,\ldots\,,\vec{x}_N)~;\label{coordinates}
\end{equation}
see Sec.~\ref{sec:BBGKY} for a discussion on the motivation for this coordinate collapse.
For Boson particles, $\psi$ remains invariant under permutations of $\xv_a$'s. Thus, the integrals
entering the sum of~(\ref{estimate-I:NB}) are the same and the estimate~(\ref{estimate-I:NB}) evidently becomes\looseness=-1
\begin{equation}
m_{ab}:=\int {\rm d}t\,{\rm d}\vec{X}_{ab}\ \big\{\rho\big(t,\vec{X}_{ab})\big\}
\leq \Vert\psi\Vert_{H^{1}}\Vert\psi\Vert_{L^{2}}~. \label{estimate}
\end{equation}

We proceed to prove~(\ref{estimate-I:NB}). Our program consists of the following steps.
(i) An evolution equation is derived for the particle momentum density,
expressing the fact that momentum is not conserved by the flow; see~(\ref{conservation of momentum}). (ii) This equation 
is contracted with suitable vector fields to yield~(\ref{contraction}), which contains ``error'' terms
depending on the nature of particle interactions. (iii) Estimates for these error terms are obtained
directly. (iv) The contracted equation is integrated in space and time to yield~(\ref{estimate-I:NB}).

We now describe the procedure in detail. The evolution equation~(\ref{evolution}) is recast to the form
\begin{equation}
i\partial_{t}\psi +\big(\sum_{a}H_{a}\big)\psi +
\big(\sum_{a\not= b}H_{ab}\big)\psi =0 \label{evolution-II}
\end{equation}
where $H_{a}=-\Delta_{a}$ and $H_{ab}=V(\vert \vec{x}_{a}-\vec{x}_{b}\vert)$,
which implies that our analysis will be independent of any particular scaling of the interaction potential.
By~(\ref{eq:repulsion}) the basic assumption is  
\begin{equation}
V^{\prime}(s)\leq 0\qquad \forall\, s\in\mathbb{R}^+~.\label{eq:repulsion-II}
\end{equation} 
The particle coordinates are denoted by $x^{j}_{a}$ where $j=1,\,2,\,3$ and $a=1,\,2,\,\ldots\,, N$.
We consider the Euclidean space with metric 
\begin{equation}
g^{ab}_{jk}:=\delta_{jk}\delta^{ab}~,\quad j,\,k=1,\,2,\,3~,\quad a,\,b=1,\,2,\,\ldots\,, N~,
\label{metric}
\end{equation}
and use the notation
\begin{equation}
\nabla^{j}_{a}:={\partial\over\partial x^{j}_{a}}~.\label{derivatives}
\end{equation}

Next, we define the momentum variables and stress tensor.
The components of the momentum density 
for the $a$th particle are defined by
\begin{equation}
p^{j}_{a}:={1\over 2i}\big(\psi^{\ast}\nabla^{j}_{a}\psi -\psi
\nabla^{j}_{a}\psi^{\ast}\big)~,\quad j=1,\,2,\,3~,\qquad a=1,\,2\ldots N~,
\label{momentum}
\end{equation}
where $\psi^\ast$ is the complex conjugate of $\psi$.
The stress tensor is
\begin{equation}
\sigma^{jk}_{ab}:=(\nabla^{j}_{a}\psi)(\nabla^{k}_{b}\psi^{\ast})
+(\nabla^{j}_{a}\psi^{\ast})(\nabla^{k}_{b}\psi)~;\quad j,\,k=1,\,2,\,3,\ a,\,b=1,\,\ldots\,,N~. 
\label{stress tensor}
\end{equation}

For completeness, we provide the conservation law for the density $\rho$,
although this law is not invoked directly in the proof. 
The mass conservation statement reads\looseness=-1
\begin{equation}
\partial_{t}\{\rho\} -\nabla^{a}_{j}\{p^{j}_{a}\} =0~. 
\label{conservation of mass}
\end{equation}
Use of this law will be made at the end of this section and in Sec.~\ref{sec:commutator}.

The evolution equation for the momenta $p^a_j$ stems from differentiation of~(\ref{momentum})
with respect to time and use of~(\ref{evolution}). The resulting equation reads
\begin{equation}
\partial_{t}\{p_{a}^{j}\}-\nabla^{k}_{b}\{\sigma^{jb}_{ka}\} +\nabla^{j}_{a}\{{\bf \Delta}\rho\}
-\sum_{b,b\not= a}2V^{\prime}(\vert \xv_{a}-\xv_{b}\vert)\ {x_{a}^j-x_{b}^j\over
\vert \xv_{a}-\xv_{b}\vert}\ \rho =0 \label{conservation of momentum}
\end{equation}
where ${\bf\Delta}:=\sum_{a}\Delta_{a}$ is the $3N$-dimensional Laplacian.
Equation (\ref{conservation of momentum}) is not a conservation law because of $V^\prime$
on the left-hand side. By defining the weights $w_{ab}$ and the $3N$-dimensional vector 
$\vec{M}_a=\big(\vec{M}^{\,1},\,\vec{M}^{\,2},\,\ldots\,,\vec{M}^{\,N}\big)$ according to\looseness=-1
\begin{equation}
w_{ab}:=-2V^{\prime}\big(\vert\vec{x}_{a}-\vec{x}_{b}\vert\big)~,\qquad 
\vec{M}_{a}:=\sum_{b,b\not= a}w_{ab}
{\vec{x}_{a}-\vec{x}_{b}\over
\vert \vec{x}_{a}-\vec{x}_{b}\vert}~, \label{weights}
\end{equation}
we rewrite~(\ref{conservation of momentum}) as
\begin{equation}
\partial_{t}\{\vec{p}^{\ a}\}-\nabla^{b}_{k}
\big\{\vec{\sigma}^{\ ka}_{b}\big\} +\vec{\nabla}^{a}\big\{{\bf \Delta}\rho\big\}
+\vec{M}^{a}\rho =0~,\qquad a=1,\,2,\,\ldots\,,N~. \label{2nd momentum conservation}
\end{equation}
For algebraic convenience we have written $\vec{p}^{\ a}:=(p^{a}_1,p^{a}_{2},p^{a}_3)$,
the three-vector with respect to the group of coordinates for $\vec{x}_{a}$. 
Notice that the condition $V^{\prime}\leq 0$ for repulsive interactions entails $w_{ab}=w_{ba}\geq 0$.

Next, we contract the momentum equation,~(\ref{conservation of momentum}) or~(\ref{2nd momentum conservation}),
with a suitable $3N$-dimensional vector field, $Y_{(ab)}$. This field consists of $N$ ordered three-vectors
and is defined by
\begin{equation}
Y_{(ab)}:=\big(\vec{0},\,\ldots\,,\vec{0},\, {\xv_{a}-\xv_{b}\over\vert \xv_{a}-\xv_{b}\vert},\, \vec{0},\,\ldots\,,
\vec{0},\, {\xv_{b}-\xv_{a}\over\vert \xv_{b}-\xv_{a}\vert},\,\vec{0},\,\ldots\,,\vec{0}\,\big)~.\label{vector fields}
\end{equation}
This expression means that the sole nonzero three-vectors forming $Y_{(ab)}$ are the ones in
the $a$th and $b$th position in the way indicated above. We adopt the
convention that $Y_{(ab)}=Y_{(ba)}$ and write $Y^{j}_{c(ab)}$ to denote the $j$th component
of the vector located in the $c$th position of $Y_{(ab)}$ ($j=1,\,2,\,3$ and $a,\,b,\,c=1,\,\ldots,\,N$). 
The contraction of the $a$th-particle momentum, $\vec{p}^{\,a}$, with $Y_{a(cd)}$ equals
\begin{equation}
p_j^aY^j_{a(cd)}=\vec{p}_{c}\cdot\frac{\xv_c-\xv_d}{|\xv_c-\xv_d|}+\vec{p}_{d}\cdot\frac{\xv_d-\xv_c}{|\xv_d-\xv_c|}~,
\label{mom-contraction}
\end{equation}
where no summation over $c$ and $d$ is implied.
The contraction of~(\ref{conservation of momentum}) with the field $Y^{j}_{a(cd)}$ by
(\ref{vector fields}) produces the equation
\begin{eqnarray}
\lefteqn{0=\partial_{t}\big\{p^{a}_{j}Y^{j}_{a(cd)}\big\}
-\nabla^{b}_{k}\big\{\sigma^{ka}_{jb}Y^{j}_{a(cd)}\big\}
+\nabla^{a}_{j}\big\{Y^{j}_{a(cd)}{\bf \Delta}\rho\big\}}\nonumber\\
&&+\big(\nabla^{b}_{k}Y^{j}_{a(cd)}\big)\sigma^{ka}_{jb}
-\big(\nabla^{a}_{j}Y^{j}_{a(cd)}\big){\bf \Delta}\rho +\sum_{c}w_{ac}{x_{j}^{a}-x_{j,c}\over
\vert \xv_{a}-\xv_{c}\vert}Y^{j}_{a(cd)}\rho. \label{contraction}
\end{eqnarray}

We will show that the dispersive nature of~(\ref{contraction}) provides an apriori bound for the 
collapsed density according to~(\ref{estimate-I:NB}). For this purpose, we integrate~(\ref{contraction})
over the time-slice $[0,T]\times \mathbb{R}^{3N}$.
The last two terms in the first line of~(\ref{contraction}) are explicitly integrated out to zero, viz.
\begin{equation}
\int {\rm d}t\,{\rm d}\vec{X}\ \nabla^{b}_{k}\big\{\sigma^{ka}_{jb}Y^{j}_{a(cd)}\big\}
=0=\int {\rm d}t\,{\rm d}\vec{X}\ \nabla^{a}_{j}\big\{Y^{j}_{a(cd)}{\bf \Delta}\rho\big\}~.\label{zero-terms}
\end{equation}
We will return to the first term, $\partial_t\big\{p^{a}_{j}Y^{j}_{a(cd)}\big\}$, near the end of this proof.

Our main task now is to investigate the nature of terms in the second line of~(\ref{contraction}). 
We sum over all contractions with the vector fields $Y_{(cd)}$ for $c\not= d$.
There are $N(N-1)/2$ different vector fields.
For simplicity we will relabel indices by $c=a$, $d=b$ and $a=c$. \looseness=-1

First, we focus on the term $-\big(\nabla^{a}_{j}Y^{j}_{a(cd)}\big){\bf \Delta}\rho$, which is integrated to furnish
\begin{equation}
\sum_{c\not= d}
\int {\rm d}t\,{\rm d}\vec{X}\ \left\{-\Big(\nabla^{a}_{j}Y^{j}_{a(cd)}\Big){\bf \Delta}\rho\right\}~,
\quad \vec{X}=(\vec{x}_{1},\,\vec{x}_{2},\,\ldots\,, \vec{x}_{N})~. 
\label{time-space-int}
\end{equation}
The identity
\begin{equation}
{\rm div} \vec{Y}_{(ab)} =\vec{\nabla}\cdot \vec{Y}_{(ab)}={4\over\vert\vec{x}_{a}-\vec{x}_{b}\vert}
\label{div-Y}
\end{equation}
and the subsequent integration by parts in~(\ref{time-space-int}) with $c=a$ and $d=b$ yield 
\begin{equation}
\sum_{a\not= b}\int {\rm d}t\,{\rm d}\vec{X}\ \left\{{-4\over\vert\vec{x}_{a}-\vec{x}_{b}\vert}\right\}
{\bf \Delta}\rho \ = c
\sum_{a\not= b}\int {\rm d}t\,{\rm d}\vec{X}_{ab}\ \big\{\rho\big(t,\vec{X}_{ab}\big)\big\}\ge 0~,
\label{int-by-parts}
\end{equation}
where ${\bf \Delta}_a(|\xv_a-\xv_b|^{-1})=4\pi\delta(\xv_a-\xv_b)$ and $\vec{X}_{ab}$ is defined in~(\ref{coordinates}).
Notice that in our derivation thus far we did not have to assume that the
particles are Bosons. However, it should be borne in mind that in the case with Bosons the results are simplified
since all these integrals are equal.

\begin{figure}
\begin{center}
\includegraphics[width=2.86in,height=3.75in]{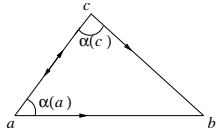}
\end{center}
\caption{Geometric interpretation of interaction terms from~(\ref{triangle1}) and~(\ref{triangle2})
by the triangle $T(\xv_a,\xv_b,\xv_c)$.}\label{F:triangle1}
\end{figure}

Next, we focus on the last term in~(\ref{contraction}), viz.
\begin{equation}
{\cal E}:=\sum_{a\not= b}
\left(\vec{M}^{a}\cdot{\vec{x}_{a}-\vec{x}_{b}\over\vert\vec{x}_{a}-\vec{x}_{b}\vert}
+\vec{M}^{b}\cdot{\vec{x}_{b}-\vec{x}_{a}\over\vert\vec{x}_{b}-\vec{x}_{a}\vert}\right)\rho(t,\vec{X})~,
\label{error}
\end{equation}
which by use of~(\ref{weights}) can be written explicitly as the triple sum
\begin{equation}
\sum_{a\not= b}\sum_{c}
\left(w_{ca}{\vec{x}_{c}-\vec{x}_{a}\over
\vert\vec{x}_{c}-\vec{x}_{a}\vert}\cdot
{\vec{x}_{b}-\vec{x}_{a}\over
\vert\vec{x}_{b}-\vec{x}_{a}\vert}+w_{cb}
{\vec{x}_{c}-\vec{x}_{b}\over
\vert\vec{x}_{c}-\vec{x}_{b}\vert}\cdot{\vec{x}_{a}-\vec{x}_{b}\over
\vert\vec{x}_{a}-\vec{x}_{b}\vert}\right)\rho(t,\vec{X})~. \label{error2-a}
\end{equation}
Terms with $c=b$ or $c=a$ in the first or second sum of~(\ref{error2-a}), respectively, add up to a manifestly positive term.
We turn our attention to terms with $c\not= a,\,b$. Let us pair the term 
\begin{equation}
w_{ca}{\vec{x}_{c}-\vec{x}_{a}\over
\vert\vec{x}_{c}-\vec{x}_{a}\vert}\cdot
{\vec{x}_{b}-\vec{x}_{a}\over
\vert\vec{x}_{b}-\vec{x}_{a}\vert}~, \label{triangle1}
\end{equation}
which comes from the contraction with $Y_{(ab)}$, with the term
\begin{equation}
w_{ac}{\vec{x}_{a}-\vec{x}_{c}\over
\vert\vec{x}_{a}-\vec{x}_{c}\vert}\cdot
{\vec{x}_{b}-\vec{x}_{c}\over
\vert\vec{x}_{b}-\vec{x}_{c}\vert}~, \label{triangle2}
\end{equation}
which comes from the contraction with $Y_{(bc)}$. The sum of these two terms
admits a geometric interpretation as follows.
Consider the triangle $T(\xv_a,\xv_b,\xv_c)$ with vertices $\vec{x}_{a}$, $\vec{x}_{b}$ and $\vec{x}_{c}$
and denote the corresponding angles $\alpha_{a}:=\alpha(a)$, 
$\alpha_{b}:=\alpha(b)$ and $\alpha_{c}:=\alpha(c)$ as shown in Fig.~\ref{F:triangle1}.
Evidently, the sum of terms~(\ref{triangle1}) and~(\ref{triangle2}) equals
\begin{equation}
w_{ac}
\big(\cos(\alpha_{a})+\cos(\alpha_{c})\big)\geq 0~,  \label{(2.20)}
\end{equation}
which has a positive sign. In a similar manner, one can pair the term
\begin{equation}
w_{cb}
{\vec{x}_{c}-\vec{x}_{b}\over\vert \vec{x}_{c}-\vec{x}_{b}\vert}
\cdot {\vec{x}_{a}-\vec{x}_{b}\over\vert\vec{x}_{a}-\vec{x}_{b}\vert}~,\label{triangle-sum1} 
\end{equation}
which comes from the contraction with $Y_{(ab)}$, with the term
\begin{equation}
w_{bc}
{\vec{x}_{b}-\vec{x}_{c}\over\vert \vec{x}_{b}-\vec{x}_{c}\vert}
\cdot
{\vec{x}_{a}-\vec{x}_{c}\over\vert\vec{x}_{a}-\vec{x}_{c}\vert} \ , \label{(2.22)}
\end{equation}
which comes from the contraction with $Y_{(ac)}$. 
The representation by the triangle of Fig.~\ref{F:triangle1} implies that the sum of the 
last two terms is
\begin{equation}
w_{cb}\big[\cos(\alpha_{b})
+\cos(\alpha_{c})\big]\geq 0~. \label{triangle-sum2}
\end{equation}

As a matter of fact, by inspection of the sum of~(\ref{triangle1}) and~(\ref{triangle2})
we can define the potential 
\begin{equation}
U\big(\vec{x}_{a},\vec{x}_{b},\vec{x}_{c}\big):=
w_{ca}{\vec{x}_{c}-\vec{x}_{a}\over
\vert\vec{x}_{c}-\vec{x}_{a}\vert}\cdot
{\vec{x}_{b}-\vec{x}_{a}\over
\vert\vec{x}_{b}-\vec{x}_{a}\vert}+
w_{ac}{\vec{x}_{a}-\vec{x}_{c}\over
\vert\vec{x}_{a}-\vec{x}_{c}\vert}\cdot
{\vec{x}_{b}-\vec{x}_{c}\over
\vert\vec{x}_{b}-\vec{x}_{c}\vert}~. \label{potential}
\end{equation}
Accordingly, the quantity ${\cal E}$ from~(\ref{error}) or~(\ref{error2-a}) is recast to the expression
\begin{equation}
{\cal E}=\sum_{a\not= b}\sum_{c}
U\big(\vec{x}_{a},\vec{x}_{b},\vec{x}_{c}\big)\rho\big(t,\vec{X}\big)~.
\label{error2}
\end{equation}

There are three remarks in order:
\vskip10pt

\noindent
{\bf Remark 2.1.} There are overall $N(N-1)^{2}$ terms in the summation 
(\ref{error2}), because each contraction
produces $2(N-1)$ terms and there are $N(N-1)/2$ contractions. The number of terms with $c=b$ or $c=a$,
which have positive sign, is $N(N-1)$; thus, the number of the remaining terms is $N(N-1)(N-2)$ and we separate them
into two groups.
\vskip5pt

\noindent
{\bf Remark 2.2.}
Notice the following term, which comes from (\ref{triangle1}):
\begin{equation}
-V^{\prime}(\vert\vec{x}_{c}-\vec{x}_{a}\vert)\ {\vec{x}_{c}-\vec{x}_{a}\over
\vert\vec{x}_{c}-\vec{x}_{a}\vert}\cdot
{\vec{x}_{b}-\vec{x}_{a}\over
\vert\vec{x}_{b}-\vec{x}_{a}\vert}\ \rho(t,\vec{X})~. \label{three-part}
\end{equation}
This term can be thought of as a three-particle interaction: the particle $c$ interacts with
$a$ and $b$. If we consider the plane $P_{a,b}$ that passes through $\vec{x}_{a}$
and is orthogonal to the vector $\vec{x}_{b}-\vec{x}_{a}$, then in the half space that
contains $\vec{x}_{b}$ the expression~(\ref{three-part}) is positive; in contrast, in the complementary
half space this expression becomes negative. The idea is that~(\ref{three-part}) can be paired with
another term where the particle $a$ interacts with $b$ and $c$ so that the sum
of the two contributions is positive. This argument gives general qualitative information
on the particle mutual repulsion.
\vskip5pt

\noindent
{\bf Remark 2.3.}  By replacing the vector fields $\vec{X}$ in~(\ref{vector fields}) with
$(\vec{x}_{a}-\vec{x}_{b}-\vec{d})/\vert\vec{x}_{a}-\vec{x}_{b}-\vec{d}\vert$, where
$\vec{d}$ is some fixed vector, we obtain an estimate like the one in~(\ref{estimate}) 
where the integration is over $\vec{x}_{a}-\vec{x}_{b}=\vec{d}$.
\vskip10pt 

These remarks conclude the investigation of the last term in evolution law~(\ref{contraction}).
We have asserted that this term is positive.

The first term in the second line of~(\ref{contraction}) reads
\begin{equation}
\Sigma :=\Big(\nabla^{b}_{k}Y^{j}_{a(cd)}\Big)\sigma^{ka}_{jb}~. \label{(2.25)}
\end{equation}
Evidently, this term is positive, i.e. $\Sigma \geq 0$ for the given vector fields. 

It remains to collect the results obtained thus far in order to derive~(\ref{estimate}).
In summary, with regard to the second line of~(\ref{contraction}) we showed that
the first and third terms are positive and the second term yields the positive integral~(\ref{int-by-parts}).
Recall that (\ref{zero-terms}) hold. By integrating the first term in~(\ref{contraction}) over
the time-slice $[0,T]\times\mathbb{R}^{3N}$ we obtain the estimate
\begin{eqnarray}
\int_0^T {\rm d}t\int{\rm d}\vec{X}\ \big\{-\partial_{t}\big(p^{a}_{j}Y^{j}_{a(cd)}\big)\big\}
&=&\int{\rm d}\vec{X}\,\big[(p_j^aY_{a(cd)}^j)\big|_{t=0}-(p_j^aY_{a(cd)}^j)\big|_{T}\big]\nonumber\\
&\le& \Vert \psi\Vert_{L^2}\,\Vert\psi\Vert_{H^1}~.\label{time-deriv-est}
\end{eqnarray}
Here, we used definition~(\ref{momentum}) for the momentum variables. This statement concludes the proof.

We close this section with the derivation of another estimate expressing the
dispersion of the particle density $\rho$. This estimate concerns the evolution of the variance of $\rho$, viz.
\begin{equation}
{{\rm d}^{2}\over {\rm d}t^{2}}\int {\rm d}\vec{X}\ \big\{ \rho D\big\}\geq 0~, \label{rho-disp}
\end{equation}
where the weight function $D(\vec{X})$ is defined by
\begin{equation}
D(\xv_{1},\,\xv_{2},\,\ldots\,, \xv_{N}):=\sum_{a\not= b}\vert\vec{x}_{a}-\vec{x}_{b}\vert~.
\label{D-def}
\end{equation}

This $D$ is intimately connected with the procedure applied above. Indeed, the derivative of $D$ reads
\begin{equation}
\vec{\nabla}^{a}D =\sum_{b\not= a}{\vec{x}_{a}-\vec{x}_{b}\over
\vert\vec{x}_{a}-\vec{x}_{b}\vert}~. \label{D-deriv}
\end{equation}
With the introduction of $Y_j^a$ by
\begin{equation}
\nabla^{a}_{j}D (\vec{x}_{1},\vec{x}_{2}\ldots \vec{x}_{N}) =: Y^{a}_{j}~,
\label{Yj-def}
\end{equation}
it becomes evident that we actually performed contractions with the vector fields $Y^a_j$ for $a=1,\,2,\,\ldots\,, N$.
Estimate~(\ref{rho-disp}) for $\rho$ results from the combination of~(\ref{D-def}) with the mass 
conservation law~(\ref{conservation of mass}).

A possible extension of our analysis is the inclusion of a nonzero trapping
potential, $V_{\tr}$. To speculate the related difficulties, we recall that the repulsive
nature of particle interactions means that the particles tend to disperse so there are no bound states.
The addition of a trapping potential in the Hamiltonian will produce an effective attractive force
that keeps the particles together. This influence will compete with the dispersive effect described above, creating the
possibility for trapped states.\looseness=-1

%
%
\section{Second Estimate by Commutator Operators}
\label{sec:commutator}
%
%

In this section we restrict attention to Bosons and derive an apriori estimate
that involves a space-time integral of the square of a reduced particle density function.
We prove this estimate by resorting to vector commutator operators and the evolution
of a suitably defined action and associated correlation function.

Specifically, we show that
\begin{equation}
\sum_{a\neq b\ }
\int\limits_{\mathbb{R}\times\mathbb{R}^{3}}{\rm d}t\, {\rm d}\vec{x}_{a-b}
\ \left\{\big(\widetilde{\rho}_{a,b}(t,\vec{x}_{ab})\big)^{2}\right\} \leq 
N^{2}\,\big\Vert\psi\big\Vert_{H^{1}}\big\Vert\psi\big\Vert^{3}_{L^{2}}~, 
\label{square-est}
\end{equation}
where the averaged, reduced density $\widetilde{\rho}_{a,b}$ is defined by
\begin{equation}
\widetilde{\rho}_{a,b}\big(t,\vec{x}_{a-b}\big)
:=\int\limits_{\mathbb{R}^{3(N-1)}} {\rm d}\vec{X}_{a,b}\, {\rm d}\vec{x}_{a+b}\ \left\{\rho\right\}~.\label{tilderho-def}
\end{equation}
In the above, the variables $\xv_{a\pm b}$ denote the center-of-mass coordinates 
\begin{equation}
\vec{x}_{a\pm b}:={1\over\sqrt{2}}\big(\vec{x}_{a}\pm \vec{x}_{b}\big)~,\label{com}
\end{equation}
and the reduced coordinates $\vec{X}_{a,b}$ are defined by
\begin{equation}
\vec{X}_{a,b}:=\big(\vec{x}_{c}\big)_{c\not= a,b}\in \mathbb{R}^{3(N-2)}~,\label{Xab-def}
\end{equation}
i.e., $\vec{X}_{a,b}$ stem from $(\xv_1,\,\xv_2,\,\ldots\,,\xv_N)$ with the pair $(\vec{x}_{a},\,\vec{x}_{b})$ being omitted. 

First, we briefly review and comment on the main assumptions and starting equations. 
The many-particle wave function satisfies~(\ref{evolution-II}), which is recast to the equation
\begin{equation}
i\partial_{t}\psi -{\bf \Delta}\psi +
\sum_{a\not= b}V\big(
\vert\vec{x}_{a}-\vec{x}_{b}\vert\big)\psi =0 \label{evolution-III}
\end{equation}
where ${\bf \Delta}$ is the Laplacian in $\mathbb{R}^{3N}$. The key idea in this section is to invoke commutator
vector operators, which we apply to conservation laws stemming from~(\ref{evolution-III}). 
In addition to condition~(\ref{eq:repulsion-II}) on $V$, we impose Bose symmetry, i.e. 
require that $\psi$ remain invariant under permutation of any coordinate pair
$(\vec{x}_{a},\vec{x}_{b})$:
\begin{equation}
\psi(t,\,\ldots\,, \vec{x}_{a},\,\ldots\,,\vec{x}_{b},\,\ldots)=
\psi(t,\,\ldots\,, \vec{x}_{b},\,\ldots\,,\vec{x}_{a},\,\ldots)~. \label{B-symmetry}
\end{equation}
This property is inherited by the density function defined by~(\ref{density}), viz.
\begin{equation}
\rho(t,\,\ldots\,, \vec{x}_{a},\,\ldots\,,\vec{x}_{b},\,\ldots)=
\rho(t,\,\ldots\,, \vec{x}_{b},\,\ldots\,,\vec{x}_{a},\,\ldots)~. 
\label{B-symmetry-rho}
\end{equation}

Crucial in our considerations are the evolution laws~(\ref{conservation of mass}) and~(\ref{conservation of momentum}) for mass 
density and momentum. These equations read
\begin{equation}
\partial_{t}\rho\ -\nabla^{b}_{k}p^{k}_{b}= 0~,\label{conserv-rho}
\end{equation}
\begin{equation}
\partial_{t}p^{a}_{j}-\nabla^{b}_{k}\left\{\sigma^{ka}_{jb}+\delta_{j}^{\ k}\delta_{b}^{\ a}
\big(-{\bf\Delta}\rho\big)\right\} +M_{j}^{a}\rho=0~,  \label{conserv-p}
\end{equation}
where $j=1,\,2,\,3$ and $a=1,\,2,\,\ldots\,, N$. The relevant quantities,
namely, the density $\rho$, the momenta $p_{j}^{a}$, and 
the stress tensor $\sigma_{jk}^{ab}$ are defined by~(\ref{density}),~(\ref{momentum})
and~(\ref{stress tensor}). In particular, $\sigma_{jk}^{ab}$ is rewritten as
\begin{equation}
\sigma_{jk}^{ab}=\rho^{-1}\big[p_{j}^{a}p_{k}^{b}+(\nabla_{j}^{a}\rho)(\nabla_{k}^{b}\rho)\big]~.
\label{stress tensor-II}
\end{equation}
It is worthwhile mentioning the conservation of the total energy, which we will not use directly here.
This law can be expressed as the conservation of the integral
\begin{equation}
E(t):=\int\limits_{\mathbb{R}^{3N}}{\rm d}\vec{X}\ \Biggl\{{1\over 2}\,{\rm tr}(\sigma)+\sum_{a\not= b}
V(\vert\vec{x}_{a}-\vec{x}_{b}\vert)2\rho\Biggr\}~.
\label{energy}
\end{equation}
This quantity is a constant (${\rm d}E/{\rm d}t\equiv 0$) and controls the right-hand side
of (\ref{square-est}) provided that $V$ is bounded below; see~(\ref{eq:repulsion}).

The proof of~(\ref{square-est}) is based on the
construction of suitable vector commutator operators and study of associated evolution laws.
Our program consists of the following steps. (i) A vector commutator operator, $\vec{C}$,
is constructed as a suitable average over all particles. 
(ii) An action, $L(t)$, is defined as the inner product $\big\langle \vec{C}\cdot \vec{P}(t)\,\big|\,\rho(t)\big\rangle$
where $\vec{P}$ is the $3N$-dimensional momentum vector. (iii) The evolution of $L(t)$
is described in terms of distinct inner products, $S_\ell$. (iv) Estimates are derived for each $S_\ell$.
(v) The evolution equation for $L(t)$ is integrated to yield~(\ref{square-est}).

We proceed to carry out this program. For each particle pair $(a,b)$, we construct the integral operator $B_{ab}$ by
\begin{equation}
\big(B_{ab}f\big)(\vec{x}_{a-b})
:=\int\limits_{\mathbb{R}^{3N}}
{{\bf 1}(\vec{x}^{\,\prime}_{a+b})\,{\bf 1}(\vec{X}^{\,\prime}_{a,b})\over
\vert\vec{x}_{a-b}-\vec{x}^{\,\prime}_{a-b}\vert}\,f\big(\vec{X}^\prime\big)\ {\rm d}\vec{X}^\prime~.
\label{Bab-def}
\end{equation}
By using these $B_{ab}$, we construct the $3N$-vector commutator operator
\begin{equation}
\vec{C}_{a;b}:=\big(0\,\ldots 0,\,\big[\vec{x}_{a-b}\,;\,B_{ab}\big],\,0\,\ldots\,0\big)~,\label{Cab-def}
\end{equation}
where $[\mathcal A\,;\,\mathcal B]:=\mathcal A\mathcal B-\mathcal B\mathcal A$ denotes the commutator
of $\mathcal A$ and $\mathcal B$, and
$[\vec{x}_{a-b}\,;\,B_{ab}]$ in~(\ref{Cab-def}) is placed in the $a$th position and acts on functions $f$ according to
\begin{equation}
\Big(\big[\vec{x}_{a-b}\,;\,B_{ab}\big]f\Big)(\vec{x}_{a-b})
=\int\limits_{\mathbb{R}^{3N}}{\vec{x}_{a-b}-\vec{x}^{\,\prime}_{a-b}\over
\vert\vec{x}_{a-b}-\vec{x}^{\,\prime}_{a-b}\vert}\ 
f\big(\vec{X}^{\prime}\big)\ {\rm d}\vec{X}^{\prime}~. \label{(3.25)}
\end{equation}
For our purposes, it is desirable to symmetrize $\vec{C}_{a;b}$ and replace it by the operator\looseness=-1
\begin{equation}
\vec{C}_{(a;b)}:=\big(0\,\ldots\, 0,\,[\vec{x}_{a-b}\,;\,B_{ab}],\,0\,\ldots\,0,\,
[\vec{x}_{b-a};B_{ba}],\,0\,\ldots\, 0\big)=\vec{C}_{a,b}+\vec{C}_{b,a}~. \label{Csymm-def}
\end{equation}
Subsequently, we average over all particles to obtain the $3N$-vector operator
\begin{equation}
\vec{C}:=\sum_{a< b}\vec{C}_{(a;b)}~. \label{Ctot-def}
\end{equation}

We use this $\vec{C}$ to construct an appropriate action $L(t)$ whose evolution
paves the way to estimate~(\ref{square-est}). With this goal in mind, let
\begin{equation}
\vec{P}:=\big\{p_{j}^{a}\big\}_{j=1,2,3}^{a=1,\ldots N} \label{TotalMomentum}
\end{equation}
be the overall momentum vector and $C^{j}_{a}$ denote the components of $\vec{C}$.
Accordingly, we consider the action
\begin{equation}
L(t):=\big\langle\vec{C}\cdot \vec{P}(t)\ \big\vert\ \rho(t)\big\rangle~. \label{L-def}
\end{equation}

The time evolution of $L(t)$ is described by
\begin{equation}
\dot{L}(t)=\big\langle\vec{C}\cdot\partial_{t}\vec{P}(t)\ \big\vert\ \rho(t)\big\rangle
-\big\langle\vec{P}(t)\ \big\vert\ \vec{C}\partial_{t}\rho(t)\big\rangle~.  \label{dot-L}
\end{equation}
By using the conservation laws~(\ref{conserv-rho}) and~(\ref{conserv-p}) and transferring the operator 
$\vec{C}{\bf \nabla}$ on the right-hand side of the inner product, we obtain the equation
\begin{equation}
\dot{L}(t)=S_{\rm cm}+S_{\rm cv}+S_{\rm ds}+S_{\rm pr}~. \label{dot-L-Sterms}
\end{equation}
The terms $S_\ell$ ($\ell=$cm, cv, ds, pr) on the right-hand side express distinct physical
effects and are defined by
\begin{eqnarray}
&&S_{\rm cm}:=\Big\langle \rho^{-1}\nabla_{j}^{a}\rho\nabla_{k}^{b}\rho
\ \big\vert\ \big(\nabla_{k}^{b}C^{j}_{a}\big)\rho\Big\rangle~, \label{Scm}\\
&&S_{\rm cv}:=\Big\langle \rho^{-1}p_{j}^{a}p^{k}_{b}\ \big\vert\
\big(\nabla_{k}^{b}C^{j}_{a}\big)\rho\Big\rangle -
\Big\langle p^{k}_{b}\ \big\vert\ \big(\nabla_{k}^{b}C^{j}_{a}\big)p_{j}^{a}\Big\rangle~,
\label{Scv} \\
&&S_{\rm ds}:=\Big\langle ({\rm div}\vec{C})(-{\bf\Delta}\rho)\ \big\vert\ \rho\Big\rangle~, \label{Sds}\\
&&S_{\rm pr}:=\big\langle\vec{C}\cdot\big(\vec{M}\rho\big)\ \big\vert\ \rho\Big\rangle~. \label{Spr}
\end{eqnarray}
The origin of these terms is described as follows.
The term $S_{\rm cm}$ is due to the compressibility of the fluid described by the
conservation laws; $S_{\rm cv}$ is a convective term; $S_{\rm ds}$ is due to dispersion; and
$S_{\rm pr}$ signifies a pressure contribution. 

Our next goal is to show that all terms $S_\ell$ in~(\ref{Scm})--(\ref{Spr}) are positive; in particular,
$S_{\rm ds}$ reduces to the integral appearing in~(\ref{square-est}). 
First, we derive alternative expressions for $S_{\rm cm}$ and $S_{\rm cv}$ in terms of appropriate tensor products;
see~(\ref{Scm-tensor}) and~(\ref{Scv-tensor}) below.
To this end, we introduce some additional formulas which will be used below.
Recall definition~(\ref{Ctot-def}) and write $\vec{C}=\big(\vec{C}_{1},\,\vec{C}_{2},\,\ldots,\,\vec{C}_{N}\big)$;
the component $\vec{C}_a$ of $\vec{C}$ reads
\begin{equation}
\vec{C}_{a}=\sum_{b, b\not= a}\big[\vec{x}_{a-b}\ ;\ B_{ab}\big] \label{eq:Ca}
\end{equation}
and has derivative
\begin{equation}
{\bf\nabla}^{c}\vec{C}_{a}=\sum_{b,b\not= a}
\Big({\bf\nabla}^{a-b}\big[\vec{x}_{a-b}\ ;\ B_{ab}\big]\Big)
\big(\delta^{c}_{\ a}-\delta^{c}_{\ b}\big)~. \label{C-deriv}
\end{equation}
The combination of the identity
\begin{equation}
\nabla_{k}^{a-b}\big[x^{j}_{a-b}\ ; B_{ab}\big]
=\delta_{k}^{\ j}B_{ab}+\big[x^{j}_{a-b}\ ;(\nabla_{k}^{a-b}B_{ab})\big] \label{B-ident}
\end{equation}
with the formula
\begin{eqnarray}
\lefteqn{\Big(\big[x^{j}_{a-b}\ ;\ \nabla_{k}^{a-b}B_{ab}\big]f\Big)(\vec{x}_{a-b})=}\nonumber\\
&&-\int\limits_{\mathbb{R}^{3N}}
{(x^{j}_{a-b}-x^{\,\prime ,j}_{a-b})(x_{k,a-b}-x^{\,\prime}_{k,a-b})\over
\vert\vec{x}_{a-b}-\vec{x}^{\,\prime}_{a-b}\vert^{3}}\,
f\big(\vec{X}^{\prime}\big)\ {\rm d}\vec{X}^{\prime}\label{B-act-f}
\end{eqnarray}
yields the expression
\begin{equation}
\Big(\nabla^{a-b}_{k}\big[x^{j}_{a-b}\ ;\ B_{ab}\big]f\Big)(\vec{x}_{a-b})
=\int\limits_{\mathbb{R}^{3N}}r_{k}^{j}\big(\vec{x}_{a-b};\vec{x}^{\,\prime}_{a-b}\big)\,
f\big(\vec{X}^{\prime}\big)\ {\rm d}\vec{X}^{\prime}~, \label{commB-f}
\end{equation}
where the kernel $r_{k}^{j}$ is defined by
\begin{eqnarray}
\lefteqn{
r_{k}^{ j}\big(\vec{x}_{a-b};\vec{x}^{\ \prime}_{a-b}\big):=
\vert\vec{x}_{a-b}-\vec{x}^{\ \prime}_{a-b}\vert^{-3} }\nonumber\\
&&\times\big[\vert\vec{x}_{a-b}-\vec{x}^{\ \prime}_{a-b}\vert^{2}\ \delta_{k}^{\ j}
-(x_{k,a-b}-x^{\prime}_{k,a-b})(x^{j}_{a-b}-x^{\prime ,j}_{a-b})\big]~.
\label{kernel-rkj-def}
\end{eqnarray}
This kernel is positive definite. By the usual tensor convention we have
\begin{equation}
{\bf r}_{ab}:=\big(r_{k}^{\ j}
(\vec{x}_{a-b};\vec{x}^{\prime}_{a-b})\big)~,\quad j,k=1,2,3~.
\label{r-tensor}
\end{equation}

With regard to the momentum vector $\vec{P}$ we write $\vec{P}=\big(\vec{p}^{\ 1},\vec{p}^{\ 2},\ldots ,\vec{p}^{\ N}\big)$
and denote the tensor product of two momentum components by
\begin{equation}
\vec{p}^{\,a}\otimes\vec{p}^{\,b}=\big(p^{a}_{j}p^{b}_{k}\big)~,\quad j,\,k=1,\,2,\,3~. 
\label{momentum-tensor-prod}
\end{equation}
For a pair $(a,b)$ of particles we form the linear combinations
\begin{equation}
\vec{p}^{\,a\pm b}:={1\over\sqrt{2}}\big(\vec{p}^{\,a}\pm \vec{p}^{\,b}\big)~,\quad a,\,b=1,\,2,\,\ldots,\, N~, 
\label{p-com}
\end{equation}
which correspond to the center-of-mass coordinates~(\ref{com}).
Consequently, we find the useful relation 
\begin{equation}
\vec{p}^{\,a-b}\otimes\vec{p}^{\,a-b}={1\over 2}\left(
\vec{p}^{\,a}\otimes\vec{p}^{\,a}+\vec{p}^{\,b}\otimes\vec{p}^{\,b}
-\vec{p}^{\,a}\otimes\vec{p}^{\,b}-\vec{p}^{\,b}\otimes\vec{p}^{\,a}\right)~. \label{p-tensor-identity}
\end{equation}
By combining~(\ref{C-deriv})--(\ref{p-tensor-identity}), we express $S_{\rm cm}$
and $S_{\rm cv}$ from~(\ref{Scm}) and~(\ref{Scv}) as
\begin{eqnarray}
S_{\rm cm}&=&\big\langle\rho^{-1}{\bf\nabla}\rho\otimes{\bf\nabla}\rho\ \big\vert\
\big({\bf\nabla}\vec{C}\big)\rho\big\rangle~, \label{Scm-tensor}\\
S_{\rm cv}&=&\big\langle\rho^{-1}\vec{P}\otimes\vec{P}\ \big\vert\
\big({\bf\nabla}\vec{C}\big)\rho\Big\rangle -
\Big\langle\vec{P}\ \big\vert\ \big({\bf\nabla}\vec{C}\big)\vec{P}\big\rangle~. \label{Scv-tensor}
\end{eqnarray}

With~(\ref{Sds}),~(\ref{Spr}),~(\ref{Scm-tensor}) and~(\ref{Scv-tensor}) in mind, we now
show that the terms $S_{\rm cm}$, $S_{\rm cv}$, $S_{\rm ds}$ and $S_{\rm pr}$ entering~(\ref{dot-L})
are all positive. In particular, the dispersive term $S_{\rm ds}$ admits a space integral representation
directly related to the left-hand side of~(\ref{square-est}).

First, we focus on the convective term, $S_{\rm cv}$. Since the tensor 
${\bf r}(\vec{x}_{a-b};\vec{x}^{\ \prime}_{a-b})$ of~(\ref{r-tensor}) is evaluated at two points, we form the associated two-point momentum
\begin{equation}
\vec{J}\big(\vec{X};\vec{X}^{\,\prime}\big)
:=\sqrt{\rho^{\,\prime}/\rho}\ \vec{P}-\sqrt{\rho/\rho^{\,\prime}}\ \vec{P}^{\,\prime}~,
\label{J-def}
\end{equation}
where
\begin{eqnarray}
&&\vec{P}:=\vec{P}\big(\vec{x}_{a+b},\,\vec{x}_{a-b},\,\vec{X}_{a,b}\big)~,\quad
\vec{P}^{\,\prime}=
\vec{P}\big(\vec{x}^{\,\prime}_{a+b},\,\vec{x}^{\,\prime}_{a-b},\,\vec{X}^{\,\prime}_{a,b}\big)~,
\label{capP-2pt} \\
&&\rho:=\rho\big(\vec{x}_{a+b},\vec{x}_{a-b},\vec{X}_{a,b}\big)~,\quad 
\rho^{\,\prime}=\rho\big(\vec{x}^{\,\prime}_{a+b},\,\vec{x}^{\,\prime}_{a-b},\,\vec{X}^{\,\prime}_{a,b}\big)~.
\label{rho-2pt} 
\end{eqnarray}
Accordingly, the convective term of~(\ref{Scv-tensor}) is expressed as the sum 
\begin{equation}
S_{\rm cv}=2\sum_{a\neq b}W_{ab}~, \label{Scv-W}
\end{equation}
where the integrals $W_{ab}$ are given by
\begin{equation}
W_{ab}:=\int
{\rm d}\vec{x}_{a+b}\,{\rm d}\vec{x}^{\,\prime}_{a+b}\,{\rm d}\vec{X}_{a,b}\,{\rm d}\vec{X}^{\,\prime}_{a,b}
\,\left\{ \Big\langle\vec{J}^{a-b}\otimes\vec{J}^{a-b}\,\big\vert\,{\bf r}_{ab}\Big\rangle_{a-b}
\right\}\geq 0 \label{Wab-def}~,
\end{equation}
with the obvious convention
\begin{equation}
\vec{J}^{\,a\pm b}:={1\over\sqrt{2}}\,\big(\vec{J}^{a}\pm\vec{J}^{b}\big)~.
\label{two point momenta} 
\end{equation}
In~(\ref{Wab-def}), $\big\langle\cdot\ \vert\ \cdot\big\rangle_{a-b}$ denotes the expectation value, or
average, over the variables $(\vec{x}_{a-b},\,\vec{x}^{\,\prime}_{a-b})$.
By~(\ref{Scv-W}) and~(\ref{Wab-def}) we conclude that 
\begin{equation}
S_{\rm cv}\geq 0~.\label{Scv-pos}
\end{equation}

Second, we concentrate on the compressible term, $S_{\rm cm}$, from~(\ref{Scm-tensor}). 
By introducing the gradient operator
\begin{equation}
{\bf\nabla}^{a\pm b}={1\over\sqrt{2}}\,\big({\bf\nabla}^{a}\pm {\bf\nabla}^{b}\big)~,
\label{grad-apmb}
\end{equation}
$S_{\rm cm}$ is recast to the form
\begin{equation}
S_{\rm cm}=2\sum_{a\not= b}
\Big\langle\rho^{-1}{\bf\nabla}^{a-b}\rho\otimes{\bf\nabla}^{a-b}\rho
\ \big\vert\ ({\bf r}_{ab})\rho\Big\rangle~.  \label{Scm-sum}
\end{equation}
Without further ado, we conclude that
\begin{equation}
S_{\rm cm}\geq 0~.
\label{Scm-pos}
\end{equation}

\begin{figure}
\begin{center}
\includegraphics[width=2.86in,height=2.75in]{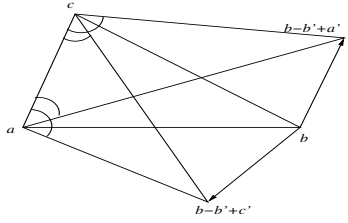}
\end{center}
\caption{Geometric interpretation of interaction terms from~(\ref{Spr-Q}) and~(\ref{Q-def}).}\label{F:triangle2}
\end{figure}

Third, we show that the pressure term, $S_{\rm pr}$, is also positive.
This case is more demanding but crucial since it is connected intimately with
the nature of particle interactions. Equation~(\ref{Spr}) is recast to the expression
\begin{equation}
S_{pr}=\sum_{a}\Big\langle\vec{C}_{a}\cdot\big(\vec{M}^{a}\rho\big)\ \big\vert\ \rho\Big\rangle
=\sum_{a,b,c\atop b\neq a,c\neq a} Q_{a,b,c}~, \label{Spr-Q}
\end{equation}
where the triple interaction terms $Q_{a,b,c}$ have the integral representation 
\begin{equation}
Q_{a,b,c}:=\int w_{ca}\,
{\vec{x}_{c}-\vec{x}_{a}\over\vert\vec{x}_{c}-\vec{x}_{a}\vert}\cdot
{\vec{x}_{b}-\vec{x}_{a}-(\vec{x}^{\,\prime}_{b}-\vec{x}^{\,\prime}_{a})
\over\vert\vec{x}_{b}-\vec{x}_{a}-(\vec{x}^{\,\prime}_{b}-\vec{x}^{\,\prime}_{a})
\vert}
\rho\big(\vec{X}\big)\,\rho\big(\vec{X}^{\,\prime}\big)\
{\rm d}\vec{X}\,{\rm d}\vec{X}^{\,\prime}~.  \label{Q-def}
\end{equation}
Recall that $w_{ca}$ is defined in~(\ref{weights}). To derive the desired estimate, we pair 
appropriately terms participating in the sum of~(\ref{Spr-Q}).
So, we pair the terms 
\begin{equation}
w_{ca}\ {\vec{x}_{c}-\vec{x}_{a}\over\vert\vec{x}_{c}-\vec{x}_{a}\vert}\cdot
{\vec{x}_{b}-\vec{x}_{a}-\vec{x}^{\,\prime}_{b}+\vec{x}^{\,\prime}_{a}
\over\vert
\vec{x}_{b}-\vec{x}_{a}-\vec{x}^{\,\prime}_{b}+\vec{x}^{\,\prime}_{a}
\vert}\ 
\rho\big(\vec{X}\big)\,\rho\big(\vec{X}^{\,\prime}\big)~, \label{pair-I}
\end{equation}
\begin{equation}
w_{ac}
{\vec{x}_{c}-\vec{x}_{a}\over\vert\vec{x}_{c}-\vec{x}_{a}\vert}\cdot
{\vec{x}_{b}-\vec{x}_{c}-\vec{x}^{\,\prime}_{b}+\vec{x}^{\, \prime}_{c}
\over
\vert
\vec{x}_{b}-\vec{x}_{c}-\vec{x}^{\,\prime}_{b}+\vec{x}^{\,\prime}_{c}
\vert}\ 
\rho\big(\vec{X}\big)\,\rho\big(\vec{X}^{\,\prime}\big)~. \label{pair-II}
\end{equation}
In the spirit of Sec.~\ref{sec:collapse}, we interpret the above terms geometrically
by considering the triangles
\begin{eqnarray}
&&T\big(\vec{x}_{a},\,\vec{x}_{c},\, \vec{x}_{b}-\vec{x}^{\,\prime}_{b}+\vec{x}^{\,\prime}_{a}
\big)=: T(a,\,c,\, b-b^{\prime}+a^{\prime})~, \label{triangle-newI} \\
&&T\big(\vec{x}_{a},\, \vec{x}_{c},\, \vec{x}_{b}-\vec{x}^{\,\prime}_{b}+\vec{x}^{\,\prime}_{c}
\big)=:T(a,\, c,\, b-b^{\prime}+c^{\prime})~,\label{triangle-newII} 
\end{eqnarray}
as shown in Fig.~\ref{F:triangle2}. By symmetry relation~(\ref{B-symmetry-rho}) for Bosons,  
we have 
$$
\rho(\ldots\, \vec{x}^{\,\prime}_{a}\,\ldots\,\vec{x}^{\,\prime}_{c}\,\ldots)
=\rho(\ldots\, \vec{x}^{\,\prime}_{c}\,\ldots\,\vec{x}^{\,\prime}_{a}\,\ldots)~.
$$
Hence, we can add to~(\ref{pair-I}) and~(\ref{pair-II}) the contributions from the terms
\begin{equation}
w_{ca}\ 
{\vec{x}_{c}-\vec{x}_{a}\over\vert\vec{x}_{c}-\vec{x}_{a}\vert}\cdot
{\vec{x}_{b}-\vec{x}_{a}-\vec{x}^{\,\prime}_{b}+\vec{x}^{\,\prime}_{c}
\over
\vert
\vec{x}_{b}-\vec{x}_{a}-\vec{x}^{\,\prime}_{b}+\vec{x}^{\,\prime}_{c}
\vert}\ 
\rho\big(\vec{X}\big)\,\rho\big(\vec{X}^{\prime \prime}\big)~,\label{w-term-I}
\end{equation}
\begin{equation}
w_{ac}\ 
{\vec{x}_{c}-\vec{x}_{a}\over\vert\vec{x}_{c}-\vec{x}_{a}\vert}\cdot
{\vec{x}_{b}-\vec{x}_{c}-\vec{x}^{\,\prime}_{b}+\vec{x}^{\,\prime}_{a}
\over
\vert
\vec{x}_{b}-\vec{x}_{c}-\vec{x}^{\,\prime}_{b}+\vec{x}^{\,\prime}_{a}
\vert}\ 
\rho\big(\vec{X}\big)\,\rho\big(\vec{X}^{\prime \prime}\big)~, \label{w-term-II}
\end{equation}
where $\vec{X}^{\prime\prime}$ denotes the variable resulting after we switch 
$\vec{x}^{\,\prime}_{a}$ and $\vec{x}^{\,\prime}_{c}$. Accordingly,
we consider again the triangles
\begin{equation}
T(a,\, c,\,b-b^{\prime}+c^{\prime})~,\quad T(a,\, c,\, b-b^{\prime}+a^{\prime})~. \label{triangles-symm}
\end{equation}
Thus, by the notation of Sec.~\ref{sec:collapse}, 
the sum of the four terms described in~(\ref{pair-I}), (\ref{pair-II}), (\ref{w-term-I})
and (\ref{w-term-II}) equals
\begin{eqnarray}
\lefteqn{w_{ca}\left\{\cos\alpha(a,c,b-b^{\prime}+a^{\prime})
+\cos\alpha(c,a,b-b^{\prime}+c^{\prime})\right\}}\nonumber\\
&&+w_{ac}\left\{\cos\alpha(a,c,b-b^{\prime}+c^{\prime})
+\cos\alpha(c,a,b-b^{\prime}+c^{\prime})\right\}~; \label{cosine-sum}
\end{eqnarray}
see Fig.~\ref{F:triangle2}.
Evidently, this term is positive. Hence, we conclude that 
\begin{equation}
S_{\rm pr}\geq 0~.\label{Spr-pos}
\end{equation}

We now focus on the dispersive term, $S_{\rm ds}$, and invoke the Fourier
transform with respect to the variables $\{\xv_a\}$. Let $\vec{\xi}_a$ be the dual variable
corresponding to the three-vector $\vec{x}_a$. By virtue of the dual $3N$-vector 
$$
\vec{\Xi}=\big(\vec{\xi}_{a+b},\, \vec{\xi}_{a-b},\, \vec{\Xi}_{a,b}\big)~,
$$
the Fourier transform of the $3N$-dimensional Laplacian reads
\begin{equation}
{\cal F}\left\{-{\bf\Delta}\right\} =
\vert\vec{\xi}_{a+b}\vert^{2}+\vert\vec{\xi}_{a-b}\vert^{2}
+\big\vert\vec{\Xi}_{a,b}\big\vert^{2}~. \label{F-Laplacian}
\end{equation}
Note that formula~(\ref{Sds}) for $S_{\rm ds}$ involves ${\rm div}\vec{C}$. 
By~(\ref{Bab-def})--(\ref{Ctot-def}), the divergence of $\vec{C}$ is
\begin{equation}
{\rm div}\vec{C} =2\sum_{a\not= b}B_{ab}~.\label{div-C}
\end{equation}
In view of~(\ref{Bab-def}), the Fourier transform of the requisite $B_{ab}$ is
\begin{equation}
{\cal F}\left\{B_{ab}\right\} ={\delta (\vec{\xi}_{a+b})\ \delta\big(\vec{\Xi}_{a,b}\big)
\over\vert\vec{\xi}_{a-b}\vert^{2}}~. \label{F-Bab}
\end{equation}
Thus, we readily compute the transform
\begin{equation}
{\cal F}\left\{{\rm div}\vec{C}(-{\bf\Delta})\right\} =\sum_{a\not= b}
{\bf 1}(\vec{\xi}_{a-b})\ \delta(\vec{\xi}_{a+b})\
\delta\big(\vec{\Xi}_{a,b}\big)~. \label{divC-Lap}
\end{equation}
By expressing $S_{\rm ds}$ in terms of Fourier transforms, 
we conclude that the dispersive term is in fact the integral
\begin{equation}
S_{\rm ds}=\sum_{a\not= b}\int\limits_{\ \mathbb{R}^{3}}
{\rm d}\vec{x}_{a-b}\ \left\{\big(\widetilde{\rho}_{ab}\big)^{2}\right\}~, \label{Sds-int}
\end{equation}
where the reduced density $\widetilde\rho\ $ is introduced in~(\ref{tilderho-def}). 

So far, we proved that the terms entering the right-hand side of~(\ref{dot-L}),
the evolution equation for $L(t)$, are positive. In particular, $S_{\rm ds}$ is described
by integral~(\ref{Sds-int}). 
The final stage of the proof involves integrating~(\ref{dot-L-Sterms}) over a time interval $[0,T]$.
Thus, we obtain the equation
\begin{equation}
\int_{0}^{T}{\rm d}t\ \left\{S_{cm}+S_{cv}+S_{ds}+S_{pr}\right\}
=L(T)-L(0)~.
\label{L-Integral}
\end{equation}
Furthermore, the action $L(t)$ defined in~(\ref{L-def}) is bounded by 
$\Vert\psi\Vert_{H^{1}}\Vert\psi\Vert_{L^{2}}^{3}$. Consequently, we reach the desired 
space-time estimate~(\ref{square-est}).
Note that the right-hand side of~(\ref{square-est}) is bounded by the energy $E(t)$ 
of the system; cf.~(\ref{energy}). This statement concludes the proof.

We close this section with a few comments. By defining the function
\begin{equation}
D\big(\vec{X};\vec{X}^{\,\prime}\big):=
\sum_{a<b}\big\vert\vec{x}_{a-b}-\vec{x}_{a-b}^{\,\prime}
\big\vert \label{Dfcn-def}
\end{equation}
and the correlation function
\begin{equation}
\mathcal C(t):=\int {\rm d}\vec{X}\,{\rm d}\vec{X}^{\,\prime}\ \left\{D\big(\vec{X};\vec{X}^{\,\prime}\big)\,
\rho(t,\vec{X})\,\rho(t,\vec{X}^{\,\prime})\right\}~,\label{cor}
\end{equation}
we deduce the relations
\begin{equation}
\dot{\mathcal C}(t)=L(t)~,\qquad \dot{L}(t)\geq 0~. \label{C-L}
\end{equation}
Because $L(t)$ is increasing and bounded by the total energy $E(t)$, we conclude that
\begin{equation}
\lim_{t\to\pm\infty}L(t) =L_{\pm}~,\label{L-limit}
\end{equation}
i.e., $\mathcal C(t)$ in~(\ref{cor}) is a convex function 
with a unique minimum, and is asymptotically linear at $t=\pm\infty$.

%
\section{Correlations and BBGKY Hierarchy and Estimate}
\label{sec:BBGKY}
%

This section has two parts. The first part serves a brief review of the description of the many-body Hamiltonian
evolution for Bosons in terms of a BBGKY hierarchy for particle reduced density matrices~\cite{esy1,esy2,esy3,esy4,spohn1,spohn2}.
In the second part, we derive an estimate for the one-particle marginal $\gamma_1$; see~(\ref{B2-estimate}) and~(\ref{fin est}) below.

\subsection{Review of BBGKY hierarchy for $k$-particle marginals}
\label{subsec:review}
 
We start again with the $N$-body Schr\"odinger equation~(\ref{evolution-III}).
The ultimate goal is to investigate aspects of this evolution for large $N$.
Following EESY~\cite{esy1,esy2,esy3,esy4} we apply the idea that the two-body interaction $V$ is scaled in some way
by incorporating the notion that pair collisions are relatively rare and weak. One particular scaling 
for $V$ is~\cite{esy1,esy2,esy3,esy4}
\begin{equation}
V(|\xv_a-\xv_b|)=N^{2}\mathcal V_1\big(N\vert \xv_{a}-\xv_{b}\vert\big) =
{1\over N}\left\{N^{3}V_1\big(N\vert \xv_{a}-\xv_{b}\vert\big)\right\}~,
\label{scl}
\end{equation}
which is equivalent to~(\ref{scattering length}) with a scattering length $l=O(N^{-1})$.
The scaling with $N$ is chosen so that the quantity inside the curly brackets approaches a constant times
the delta function as $N\to \infty$. We remind the reader that because of Boson symmetry
the $N$-body wave function $\psi$ is invariant under permutations of the space variables, i.e.
 \begin{equation}
\psi(t,\xv_1,\,\xv_2,\,\ldots\,,\xv_N)=\psi(t,\xv_{\pi(1)},\,\xv_{\pi(2)},\,\ldots\,,\xv_{\pi(N)})~,\label{Boson-s}
\end{equation}
where $\pi(a)$ is a permutation of the indices $\{1,2,\,\ldots\,,N\}$. Although we are interested in the
time evolution of the system, for the sake of convenience we will omit the time ($t$) dependence whenever this is not
relevant to the discussion.

A description of the particle system via a BBGKY hierarchy was proposed by Spohn~\cite{spohn1,spohn2} and forms the
starting point in the analysis by EESY~\cite{esy1,esy2,esy3,esy4}. A key element of this methodology
is the density matrix $\gamma_N:=|\psi\rangle\langle\psi|$, which is a trace class operator represented by
\begin{equation}
\gamma_{N}\big(t, \xv_{1},\ldots \xv_{N}\ \vert\ \xv^{\,\prime}_{1}\ldots \xv^{\,\prime}_{N}\big)
:=\psi(t, \xv_{1},\ldots \xv_{N})\ \psi^{\ast}(t, \xv^{\,\prime}_{1}\ldots \xv^{\,\prime}_{N})~.\label{gammaN}
\end{equation}
By virtue of~(\ref{evolution-III}), the evolution equation for $\gamma_N$ is
\begin{equation}
i\partial_t\gamma_N+\big(-{\bf \Delta}_N+{\bf \Delta}_{N'})\gamma_N+\sum_{a\neq b,\,a'\neq b'}(V_{ab}-V_{a'b'})\gamma_N=0~,
\label{gammaN-evolution}
\end{equation}
where ${\bf \Delta}_N=\sum_a\Delta_a$ denotes the $3N$-dimensional Laplacian, $V_{ab}:=V(|\xv_a-\xv_b|)$,
and the primed indices indicate the coordinates $(\xv^{\,\prime}_{1}\ldots \xv^{\,\prime}_{N})$.

Next, we describe features of the $k$-particle marginals $\gamma_k$ stemming from $\gamma_N$.
For this purpose, it is convenient to define $\vec{X}_{N}:=(\xv_{1},\,\xv_{2},\,\ldots\,,\xv_{N})$ and
$\vec{X}^{\,\prime}_{N}:=(\xv^{\,\prime}_{1},\,\ldots,\,\xv^{\,\prime}_{N})$.
In order to average over some of the coordinates, we also define~\cite{esy1,esy2,esy3,esy4}
\begin{equation}
\vec{X}_{k}=(\xv_{1},\,\ldots\,, \xv_{k})~,\quad \vec{X}_{N-k}=(\xv_{k+1},\,\ldots\,, \xv_{N})~; \label{CapX-def}
\end{equation}
thus, $\vec{X}_N=(\vec{X}_k,\,\vec{X}_{N-k})$ for $k=1,\,2,\,\ldots\,,N$.
The $k$-particle marginals as formed via the partial averaging 
\begin{equation}
\gamma_{k}\big(\vec{X}_{k}\ \big\vert\ \vec{X}^{\,\prime}_{k}\big) =
\int {\rm d}\vec{Y}_{N-k}\ \big\{\gamma_{N}\big(\vec{X}_{k},\vec{Y}_{N-k}\ \big\vert\ \vec{X}^{\,\prime}_{k},\vec{Y}_{N-k}\big)
\big\}~. \label{gammak-def}
\end{equation}
Note that because of~(\ref{Boson-s}) it does not matter which variables we average out. 

We now describe the evolution law for $\gamma_k$ using~(\ref{gammaN-evolution}) and definition~(\ref{gammak-def}). 
By~(\ref{com}), let $\nabla_{a\pm b}$ be the
gradient operator corresponding to $\xv_{a\pm b}$. In view of 
$\Delta_{a}-\Delta_{a^{\prime}}=\nabla_{a+a^{\prime}}\cdot\nabla_{a-a^{\prime}}$,
setting $a=a^{\prime}$ and integrating over $a+a^{\prime}$ for the averaged variables $\gamma_k$ yields zero.
Define the sets $J_{k}:=\{1,2,\,\ldots\,, k\}$,  $J^{\prime}_{k}:=\{1^{\prime},\,2^{\prime},\,\ldots\,, k^{\prime}\}$ and 
$J^{c}_{k}=\{k+1,\ldots N\}$, the complement of $J_{k}$. 
The above averaging also produces zero for the
potentials $V_{ab}$ and $V_{a'b'}$ if $a,b\in J^{c}_{k}$ or $a^{\prime},b^{\prime}\in J^{c\,\prime}_{k}$ when we identify the
variables $a=a^{\prime}$ and $b=b^{\prime}$, i.e. we identify the coordinates
$(\xv_{k+1},\,\ldots\,, \xv_{N})$ with $(\xv^{\,\prime}_{k+1},\,\ldots\,, \xv^{\,\prime}_{N})$. 
Consequently, the evolution equation for $\gamma_k$ reads
\begin{eqnarray}
\lefteqn{i\partial_{t}\gamma_{k}-\big({\bf \Delta_{k}}-{\bf \Delta}^{\prime}_{k}\big)\gamma_{k}
+\left(\sum_{a,b\in J_k\atop a\not= b}
V_{ab}-\sum_{a^{\prime}, b^{\prime}\in J^{\prime}_{k}\atop a^\prime\neq b^\prime}V_{a^{\prime}b^{\prime}}
\right) \gamma_{k}}\nonumber\\
&&+\sum_{a\in J_{k},a^{\prime}\in J^{\prime}_{k}}\sum_{\ b\in J^{c}_{k} }
\int {\rm d}\vec{Y}_{N-k}\ 
\big[V(\vert \xv_{a}-\vec{y}_{b}\vert)-V(\vert \xv_{a^{\prime}}-\vec{y}_{b}\vert)\big]\nonumber\\
&&\times\gamma_{N}\big(\xv_{1}\ldots \xv_{a}\ldots\xv_{k},\vec{y}_{k+1}\ldots \vec{y}_{b}\ldots\vec{y}_{N}\big\vert
\xv^{\,\prime}_{1}\ldots\xv^{\,\prime}_{a}\ldots\xv^{\,\prime}_{k},\vec{y}_{k+1}\ldots \vec{y}_{b}\ldots 
\vec{y}_{N}\big)\nonumber\\
\qquad&&\qquad =0~.\label{gammak-evolution} 
\end{eqnarray}

We proceed to write~(\ref{gammak-evolution}) as a BBGKY-type hierarchy for $\gamma_k$.
Because of Boson symmetry, integrations over $b\in J^{c}_{k}$ are reduced to
\begin{eqnarray}
&&(N-k)\sum_{a\in J_{k}}
\int {\rm d}\vec{Y}_{N-k-1}\,{\rm d}\vec{y}_{k+1}\ \big[V(\vert \xv_{a}-\vec{y}_{k+1}\vert)-V(\vert \xv^{\,\prime}_{a}-\vec{y}_{k+1}\vert)\big]\nonumber\\
&&\times\gamma_{N}\big(\xv_{1}\ldots \xv_{a}\ldots \xv_{k},\vec{y}_{k+1},\vec{Y}_{N-k-1}\big\vert
\xv^{\,\prime}_{1}\ldots \xv^{\,\prime}_{a}\ldots \xv^{\,\prime}_{k},\vec{y}_{k+1},\vec{Y}_{N-k-1}\big).
\label{gammak-integ} 
\end{eqnarray}
With the definitions $H_{k}:=-{\bf \Delta}_{k}+\sum_{a, b\in J_{k}}V_{ab}$ where $a\neq b$ and
\begin{eqnarray}
{\cal C}^{a,a^{\prime}}_{V,k+1}\big[\gamma_{k+1}\big]&=&
\int {\rm d}\vec{y}_{k+1}\ \big[V(\vert \xv_{a}-\vec{y}_{k+1}\vert)
-V(\vert \xv^{\,\prime}_{a}-\vec{y}_{k+1}\vert)]\big)\nonumber\\
&&\times\gamma_{k+1}\big(X_{k},y_{k+1}\vert X^{\prime}_{k},y_{k+1}\big) \ , \label{(1.16)}
\end{eqnarray}
(\ref{gammak-evolution}) takes the form 
\begin{equation}
i\partial_{t}\gamma_{k}+\big(H_{k}-H_{k^{\prime}}\big)\gamma_{k} +
(N-k)\sum_{a\in J_{k}}\sum_{\,a^{\prime}\in J^{\prime}_{k}}
{\cal C}^{a,a^{\prime}}_{V,k+1}\big[\gamma_{k+1}\big]=0 ~. \label{hierarchy}
\end{equation}
This set of coupled of equations constitutes a finite BBGKY-type hierarchy for $\gamma_k$ where
$k=1,\,2,\,\ldots\,,N$. Evidently, each $\gamma_k$ depends on $N$.

\subsection{Estimate for one-particle density matrix $\gamma_1$}
\label{subsec:1p-estimate}

Thus far, we have essentially reviewed the formulation by EESY~\cite{esy1,esy2,esy3,esy4}. In the
following, we exploit this framework to derive an estimate that involves $\gamma_1$. We assume
that the interaction potential $V_1$ of~(\ref{scattering length}) is integrable, $V_1\in L^1(\mathbb{R}^3)$. \looseness=-1

For finite $N$, the first equation of hierarchy~(\ref{hierarchy}) reads
\begin{equation}
i\partial_{t}\gamma_{1}-\big(\Delta_{1}-\Delta_{1^{\prime}}\big)\gamma_{1}
+(N-1){\cal C}^{1,1^{\prime}}_{V,2}\big[\gamma_{2}\big]=0~,\label{gamma1-eq}
\end{equation}
where $\gamma_2=\gamma_{2}\big(\xv_{1},\xv_{2}\ \vert\ \xv^{\,\prime}_{1},\xv^{\,\prime}_{2}\big)$.
To motivate the analysis for finite $N$ given below, we first consider the simpler case where $N\to\infty$, in which
the interaction potential approaches a constant times a delta function.
As $N\to\infty$, by the scaling of~(\ref{scl}) we have the limit
\begin{equation}
V_N(\cdot):=(N-1)N^{2}V_1\big(N\vert \xv_{a}-\vec{y}_{k+1}\vert\big)\to g\delta(\vert \xv_{a}-\vec{y}_{k+1}\vert)~, \label{scaling-limit}
\end{equation}
where $g$ is a fixed positive parameter ($g> 0$). In effect,
we arrive at the basic collapse mechanism that reduces $\gamma_{2}$ to a function
of $6$ variables, viz.
\begin{equation}
{\cal B}_{2}\big[\gamma_{2}\big]:=g\left\{
\gamma_{2}\big(\xv_{1},\xv_{1}\ \vert\ \xv^{\,\prime}_{1},\xv_{1}\big)
-\gamma_{2}\big(\xv_{1},\xv^{\,\prime}_{1}\ \vert\ \xv^{\,\prime}_{1},\xv^{\,\prime}_{1}\big)
\right\}~. \label{B2-collapse}
\end{equation}
Thus,~(\ref{gamma1-eq}) becomes
\begin{equation}
i\partial_{t}\gamma_{1}-\big(\Delta_{1}-\Delta_{1^{\prime}}\big)\gamma_{1} +
{\cal B}_{2}\big[\gamma_{2}\big]=0~. \label{gamma1-collapse}
\end{equation}

The partial Fourier transform of~(\ref{gamma1-collapse}) has the relatively simple structure of a transport law; see~(\ref{transport}).
Next, we derive an estimate for the forcing term ${\cal B}_2[\gamma_2]$ of~(\ref{gamma1-collapse}). Notice that the operator 
${\cal B}_{2}$ acting on $\gamma_{2}$ is written as
\begin{eqnarray}
{\cal B}_{2}\big[\gamma_{2}\big]&=&g\int {\rm d}\vec{Y}_{N-2}
\,\left\{\psi\big(\xv_{1},\xv_{1};\vec{Y}_{N-2}\big)\,\psi^{\ast}\big(\xv_{1},\xv^{\,\prime}_{1};\vec{Y}_{N-2}\big)
-\psi\big(\xv_{1},\xv^{\,\prime}_{1};\vec{Y}_{N-2}\big)\right.\nonumber\\
&&\qquad\qquad\times\left.\psi^{\ast}\big(\xv^{\,\prime}_{1},\xv^{\,\prime}_{1};\vec{Y}_{N-2}\big)\right\}~.
 \label{B2-collapse-II}
\end{eqnarray}
By writing
\begin{equation}
\widetilde\rho(\xv_{1},\xv_{2}):=\int {\rm d}\vec{Y}_{N-2}\ \left\{{1\over 2}
\big\vert\psi\big(\xv_{1},\xv_{2},\vec{Y}_{N-2}\big)\big\vert^{2}\right\} \label{red-rho-def}
\end{equation}
and denoting ${\cal B}_2[\gamma_2]$ by ${\cal B}_{2}(\xv_{1},\xv^{\,\prime}_{1})$, we have the estimate
\begin{eqnarray}
\Big\Vert{\cal B}_{2}(t)\Big\Vert^2_{L^{2}_{1-1^{\prime}}(L^{1}_{1+1^{\prime}})}&=&
\int {\rm d}\xv_{1-1^{\prime}}\ \left[
\int {\rm d}\xv_{1+1^{\prime}}
\Big\vert{\cal B}_{2}(\xv_{1+1^{\prime}},\xv_{1-1^{\prime}})\Big\vert
\right]^{2}\nonumber\\
&\leq& 2\int {\rm d}\xv_{1-1^{\prime}}
\Big[\int {\rm d}\xv_{1+1^{\prime}}\Big\{\sqrt{\widetilde\rho}\big(\xv_{1},\xv_{1}\big)
\sqrt{\widetilde\rho}\big(\xv_{1},\xv^{\,\prime}_{1}\big)\Big\}\Big]^{2}\nonumber\\
&\leq&
\left(\int {\rm d}\xv_{1}\ \widetilde\rho(\xv_{1},\xv_{1})\right)
\left[\int {\rm d}\xv_{1-1^{\prime}}\right.\nonumber\\
\mbox{}&&\times\left.\Big(\int {\rm d}\xv_{1+1^{\prime}}\ \widetilde\rho(\xv_{1},\xv^{\,\prime}_{1})\Big)\right]\nonumber\\
&=& c\int{\rm d}\xv_1\ \widetilde\rho(\xv_1,\xv_1)~, \label{B2-estimate}
\end{eqnarray}
by the definition of $\rho$ and $\widetilde\rho$.
It should be borne in mind that the density $\widetilde\rho$ depends on $t$.

In order to study the implication of~(\ref{B2-estimate}), let us write 
$\Delta_{1}-\Delta_{1^{\prime}}=\nabla_{1+1^{\prime}}\cdot\nabla_{1-1^{\prime}}$
and ${\cal B}_{2}(\xv_{1},\xv_{1^{\,\prime}})={\cal B}_{2}(\xv_{1+1^{\prime}},\xv_{1-1^{\prime}})$.
We proceed to derive an estimate involving $\gamma_1$.
By taking the Fourier transform of~(\ref{gamma1-collapse}) in the $\xv_{1+1^{\prime}}$ variable and denoting 
the dual variable by $\vec{v}$, we arrive at the transport equation
\begin{equation}
\partial_{t}\widehat{\gamma}_{1}-2\pi \vec{v}\cdot\nabla_{1-1^{\prime}}\widehat{\gamma}_{1}
=\widehat{\cal B}_{2}\big(t,\vec{v},\xv_{1-1^{\prime}}\big)~, \label{transport}
\end{equation}
where $\widehat{\gamma}_{1}(t,\vec{v},\xv_{1-1^{\prime}})$ denotes the Fourier-transformed function.
For the sake of some simplicity, we assume that the initial data for $\gamma_1$ are zero, i.e. $\gamma_{1}(t=0)=0$.
Thus, the solution to~(\ref{transport}) is
\begin{equation}
\widehat{\gamma}_{1}(t,\vec{v},\xv_{1-1^{\prime}})=
\int_{0}^{t}\widehat{\cal B}_{2}\big(s,\vec{v},\xv_{1-1^{\prime}}-2\pi s\vec{v}\big)\,{\rm d}s~. 
\label{gamma1-solution}
\end{equation}
On the basis of~(\ref{gamma1-solution}), we derive the estimate
\begin{eqnarray}
\sup_{t,\vec{v}}\Vert\widehat{\gamma}_{1}\Vert^2_{L^{2}}
&=&\sup_{t,\vec{v}}\Big\Vert\int_{0}^{t}{\rm d}s
\int {\rm d}\vec{y}\ \left\{e^{i2\pi \vec{y}\cdot\vec{v}}{\cal B}_{2}\big(s,\vec{y},\xv_{1-1^{\prime}}-2\pi s\vec{v}\big)\right\}
\Big\Vert^2_{L^{2}_{1-1^{\prime}}}\nonumber\\
&\leq& \int_0^t {\rm d}s\ \Big\Vert{\cal B}_{2}(s)\Big\Vert^2_{L^{2}_{1-1'}(L^{1}_{1+1'})}~.
\label{fin est}
\end{eqnarray}
The combination of~(\ref{fin est}) with~(\ref{B2-estimate}) and~(\ref{estimate-I:NB})
provides the desired estimate for $\gamma_1$ when $N\to\infty$.

This analysis applies directly to the case of finite $N$ with minor modifications.
First, (\ref{transport})--(\ref{fin est}) remain intact under the replacement $(N-1){\cal C}^{1,1'}_{V,2}\equiv {\cal B}_2$
in~(\ref{gamma1-eq}).
Second, in view of Remark~2.3 of Sec.~\ref{sec:collapse} and for $V_1\in L^1(\mathbb{R}^3)$, estimate~(\ref{B2-estimate}) becomes
\begin{eqnarray}
\Big\Vert{\cal B}_{2}(t)\Big\Vert^2_{L^{2}_{1-1^{\prime}}(L^{1}_{1+1^{\prime}})}
&\leq&
\left(\int{\rm d}\vec{d}\, V_N(\vec{d})\int {\rm d}\xv_{1}\ \widetilde\rho(\xv_{1},\xv_{1}+\vec{d})\right)
\left[\int {\rm d}\xv_{1-1^{\prime}}\right.\nonumber\\
\mbox{}&&\times\left.\Big(\int {\rm d}\xv_{1+1^{\prime}}\ \widetilde\rho(\xv_{1},\xv^{\,\prime}_{1})\Big)\right]\nonumber\\
&=& c\int{\rm d}\vec{d}\,V_N(\vec{d})\int {\rm d}\xv_1\ \widetilde\rho(\xv_1,\xv_1+\vec{d})\nonumber\\
&\leq& \Vert V_N\Vert_{L^1}\,\sup_{\vec{d}}\int{\rm d}\xv_1\,\widetilde\rho(\xv_1,\xv_1+\vec{d})~.\label{B2-estimate-mod}
\end{eqnarray}

It is natural to ask whether this analysis can be extended to $k$-particle marginals
for $k>1$. However, evolution equations~(\ref{hierarchy}) satisfied by $\gamma_{k}$ cannot be converted to simple transport laws if $k>1$ because of the
presence of interaction potentials. Hence, the procedure applied hitherto is not directly applicable
to the case with $k> 1$. The derivation of estimates for $\gamma_{k>1}$ is left for future work.

\end{document}